\documentclass[10pt,a4paper]{article}
\usepackage[utf8]{inputenc}
\usepackage[english]{babel}
\usepackage{cite}
\usepackage{csquotes}
\usepackage{amsmath}
\usepackage{amsfonts}
\usepackage{amssymb}
\usepackage{amsthm}
\newtheorem{theorem}{Theorem}[subsection]
\newtheorem{lemma}[theorem]{Lemma}
\newtheorem{definition}[theorem]{Definition}

\newtheorem{axiom}[theorem]{Axiom}
\usepackage{hyperref}
\hypersetup{
    colorlinks,
    citecolor=blue,
    filecolor=blue,
    linkcolor=blue,
    urlcolor=blue
}
\author{Adrian Fellhauer}
\title{On the relation of three theorems of analysis to the axiom of choice}
\begin{document}

\maketitle

\begin{abstract}
In what follows, essentially two things will be accomplished: First, it will be proven that a version of the Arzel\`a--Ascoli theorem and the Fr\'echet--Kolmogorov theorem are equivalent to the axiom of countable choice for subsets of real numbers. Secondly, some progress is made towards determining the amount of axioms that have to be added to the Zermelo--Fraenkel system so that the uniform boundedness principle holds.
\end{abstract}

\tableofcontents

\section{Introduction}

\subsection{Outline}

In subsection \ref{subsec:context_and_motivation}, we briefly describe the greater goal of this paper.

In section \ref{sec:arzela-ascoli_and_frechet-kolmogorov}, we treat the theorems of Arzel\`a--Ascoli and Fr\'echet--Kolmogorov in the following manner: In subsection \ref{subsec:arzela-ascoli_and_frechet-kolmogorov_introduction}, we state precisely the theorems which we will investigate and briefly comment on them, in subsection \ref{subsec:arzela-ascoli_and_frechet-kolmogorov_modified_versions} we prove modified versions of both theorems that hold without any choice axiom, and \emph{using these}, in subsection \ref{subsec:arzela-ascoli_and_frechet-kolmogorov_relation_to_ac} we prove that both theorems under consideration are equivalent to the axiom of countable choice for subsets of real numbers.

In section \ref{sec:ubp}, we investigate the uniform boundedness principle as follows: In subsection \ref{subsec:ubp_introduction} we state the uniform boundedness principle, clarify how it relates (in ZF) to several other theorems and state the knowledge regarding its choice-axiomatic strength until this paper, in subsection \ref{subsec:ubp_implied_by_cc} we give a new proof of the fact that the uniform boundedness principle follows from the axiom of countable choice, in subsection \ref{subsec:ubp_implied_axioms} we deduce several choice-like axioms from the uniform boundedness principle, in subsection \ref{subsec:weak_ubp_equiv_cmc} we give a weak version of the uniform boundedness principle which is equivalent to the axiom of countable multiple choice, and in subsection \ref{subsec:ubp_outlook} we elaborate on how our results are incomplete and what seem promising directions for further investigation.

<\subsection{Context and motivation}\label{subsec:context_and_motivation}

Many mathematicians blindly accept the axiom of choice. This may be because it is a straightforward generalisation of things which hold trivially in the finite to the infinite. However, history has proven (for instance in the case of summation) that such generalisations may yield contradictions. Now for the axiom of choice, the situation is somewhat different because it is \emph{logically independent} of the Zermelo--Fraenkel system (or for short: ZF); from this follows that ZF plus the axiom of choice can only lead to a contradiction if ZF already leads to a contradiction. However, if physical reality is used as a model for the axioms that are used, the Banach--Tarski paradox doubtlessly contradicts the preservation of energy and mass (since the amount of energy used in reassembling an object of a certain sufficient size is surely far lower than the energy needed for creating an entire object of the same size).

If the axiom of choice is not accepted, an alternative approach may be, instead of proving theorems for \emph{all}, say, Banach spaces or rings or whatever object, to only prove these theorems for classes of spaces which are defined such that certain choice-axiomatic properties hold for them. For instance, the proof of theorem \ref{cc_implies_ubp} and the ensuing remark will demonstrate that one may prove the uniform boundedness principle in ZF for all Banach spaces in which every countable product of open subsets is nonempty. But it is most certainly sufficient to assume that \emph{all} countable products of \emph{any} sets are nonempty (the axiom of countable choice).

In this context, the aim of this paper is to advance knowledge on exactly how much choice is needed to hold for a given space so that certain theorems are true.

\section{The Arzel\`a--Ascoli and Fr\'echet--Kolmogorov theorems}
\label{sec:arzela-ascoli_and_frechet-kolmogorov}

\subsection{Introduction}
\label{subsec:arzela-ascoli_and_frechet-kolmogorov_introduction}

When $X$ is a topological space and $Y$ is a metric space with metric $d_Y$, the set of continuous, bounded functions from $X$ to $Y$ with metric
\begin{equation*}
d(f, g) := \sup_{x \in X} d_Y(f(x),g(x)).
\end{equation*}
is a metric space which is denoted by $\mathcal C(X,Y)$. We write $\mathcal C(X)$ for $\mathcal C(X, \mathbb R)$.

There are several closely related results that bear the name \enquote*{Arzel\`a--Ascoli theorem} (see for instance Yosida \cite[p.~85f.]{yosida}, Tao \cite[Theorem~1.8.23, p.~114]{taoReal}, Brezis \cite[Theorem~4.25, p.~111]{brezis} or Rudin \cite[Theorem~7.25, p.~158]{rudin}); all concern (relative) compactness in spaces of continuous functions. We shall be concerned with the following version:

\begin{theorem}[Arzel\`a--Ascoli]\label{arzela_ascoli}
Let $K \subset \mathbb R^d$ be compact (bearing the topology on $K$ that is induced by the Euclidean topology of $\mathbb R^d$) and let $\mathcal F \subseteq \mathcal C(K)$. Then the following two are equivalent:
\begin{enumerate}
\item Every sequence in $\mathcal F$ contains a convergent subsequence.
\item $\mathcal F$ is uniformly bounded and equicontinuous.
\end{enumerate}
\end{theorem}

We will prove that theorem \ref{arzela_ascoli} is equivalent to the axiom of countable choice for subsets of the real numbers. (Note that other versions of the Arzel\`a--Ascoli theorem have also been studied with regard to their axiomatic strength, e.g. in Herrlich \cite{herrlichAscoli}.)

The Fr\'echet--Kolmogorov theorem concerns (relative) compactness in certain $L^p$ spaces; it is contained within several (perhaps most) introductory functional analysis textbooks (for instance in Brezis \cite[Theorem~4.26, p.~111]{brezis} or Yosida \cite[p.~275]{yosida}).

\begin{theorem}[Fr\'echet--Kolmogorov]
Let $1 \le p < \infty$, let $S \subset \mathbb R^d$ be bounded and measurable and let $\mathcal F \subseteq L^p(S)$. Then the following are equivalent:
\begin{enumerate}
\item Every sequence in $\mathcal F$ contains a convergent subsequence.
\item $\mathcal F$ is bounded in $L^p(S)$ and \begin{equation*}
\lim_{h \to 0} \sup_{f \in \mathcal F} \int_S |f(x +h) - f(x)|^p dx = 0.
\end{equation*}
\end{enumerate}
\end{theorem}

Note that for the integral in the above limit to make sense, $f \in \mathcal F$ is extended to $\mathbb R^d$ by being zero outside $S$.

\subsection{Modified versions of both theorems}
\label{subsec:arzela-ascoli_and_frechet-kolmogorov_modified_versions}

As a first step of investigating the choice-axiomatic nature of the Arzel\`a--Ascoli and Fr\'echet--Kolmogorov theorems, we establish modified versions of both theorems which hold true in ZF, without assuming any version of the axiom of choice. Oddly enough, these modified versions will later be needed in determining the choice-axiomatic strength of the full theorems of Arzel\`a--Ascoli and Fr\'echet--Kolmogorov as given above.

\begin{theorem}[Modified Arzel\`a--Ascoli theorem]\label{modified_ascoli}
Let $X$ be a compact\footnote{by \enquote*{compact} we mean \enquote*{every open cover has a finite subcover.} Note that with this definition the statement \enquote*{A pseudometric space is compact if and only if it is sequentially compact} is equivalent to the axiom of countable choice (see Bentley and Herrlich \cite[Theorem~4.3, p.~161]{bentleyHerrlich}).} and separable\footnote{Note that the assertion \enquote*{every compact pseudometric space is separable} is equivalent to the axiom of countable choice (see Bentley and Herrlich \cite[Theorem~4.11, p.~164]{bentleyHerrlich}).} metric space with metric $d_X$, let $Y$ be a metric space with metric $d_Y$ and let $\mathcal F \subseteq \mathcal C(X)$. Then the following two are equivalent:
\begin{enumerate}
\item Every sequence in $\mathcal F$ contains a convergent subsequence.
\item Every countable subset of $\mathcal F$ is pointwise relatively sequentially compact and equicontinuous.
\end{enumerate}
\end{theorem}

\begin{proof}
For $(2) \Rightarrow (1)$, the proof given in Yosida \cite[p.~85f.]{yosida} is essentially sufficient; note only that
\begin{itemize}
\item the Bolzano--Weierstra{\ss} theorem may be proven constructively (see for instance \cite[Theorem~2.5.5, p.~64]{abbott}), and
\item when $\{s_n | n \in \mathbb N\}$ is a countable, dense subset of $X$, then for every $\epsilon > 0$ there must automatically exist (a minimal) $k(\epsilon)$ such that
\begin{equation*}
\sup_{s \in S} \inf_{1 \le j \le k(\epsilon)} d(s_j, s) \le \epsilon;
\end{equation*}
this follows by considering the open cover $\{B_\epsilon(s_n)|n \in \mathbb N \}$ and applying compactness.
\end{itemize}

For $(1) \Rightarrow (2)$, one may consult the proof given in Tao \cite[Proof of theorem~1.8.23, $(i) \Rightarrow (ii)$, p.~114]{taoReal}; note only that in a separable space, a countable dense subset $\{x_n|n \in \mathbb N\}$ yields a choice function on the set of all open sets, for one may take the first $x_n$ contained within a given open set.
\end{proof}

\begin{theorem}[Modified Fr\'echet--Kolmogorov theorem]\label{modified_frechet}
Let $1 \le p < \infty$, let $S \subset \mathbb R^d$ be bounded and measurable and let $\mathcal F \subseteq L^p(S)$. Then the following are equivalent:
\begin{enumerate}
\item Every sequence in $\mathcal F$ contains a convergent subsequence.
\item Every countable $\mathcal G \subseteq \mathcal F$ is bounded in $L^p(S)$ and satisfies \begin{equation*}
\lim_{h \to 0} \sup_{f \in \mathcal G} \int_S |f(x +h) - f(x)|^p dx = 0.
\end{equation*}
\end{enumerate}
\end{theorem}

Note again that $f$ is extended to $\mathbb R^d$ by being zero outside $S$.

\begin{proof}
For $(1) \Rightarrow (2)$, the countability of $\mathcal G$ makes certain that the argument given in Yosida \cite[p.~275f.]{yosida} essentially goes through in ZF.

For $(2) \Rightarrow (1)$, we use our modified version of the Arzel\`a--Ascoli theorem (theorem \ref{modified_ascoli}). To do so, we assimilate elements of the standard proof given for instance in Brezis \cite[Proof of theorem~4.26, p.~111ff.]{brezis}, but transform the argument to an argument of more \enquote*{sequential} flavour. In our argument, we will use the functions $\rho_n$ that for $n \in \mathbb N$ are given by
\begin{equation*}
\rho_n : \mathbb R^d \to \mathbb R, ~~ \rho_n(x) := \begin{cases}
\frac{\displaystyle{n^d e^{\frac{1}{1 - \|nx\|^2}}}}{\displaystyle\int_{B_1(0)} e^{\frac{1}{1 - \|y\|^2}} dy } & \|x\| < 1 \\
0 & \text{else}
\end{cases}.
\end{equation*}
These functions satisfy $\rho_n \in \mathcal C^\infty(\mathbb R^d)$, $\operatorname{supp} \rho_n \subseteq \overline{B_{1/n}(0)}$ and $\int_{B_{1/n}(0)} \rho_n(x) dx = 1$ (see for instance Brezis \cite[p.~108]{brezis}).

Now let $(f_n)_{n \in \mathbb N}$ be a sequence in $\mathcal F$. Then $\mathcal G = \{f_n | n \in \mathbb N\}$ is a countable subset of $\mathcal F$. For $n, m \in \mathbb N$ define
\begin{equation*}
h_{n,m} := f_n * \rho_m.
\end{equation*}
Then define for $m \in \mathbb N$
\begin{equation*}
\mathcal H_m := \{h_{n,m} |n \in \mathbb N\}.
\end{equation*}
We claim that for each $m \in \mathbb N$ we have that $\mathcal H_m$ is uniformly bounded and equicontinuous. Indeed, uniform boundedness follows from
\begin{equation*}
\left| h_{n,m}(x) \right| \le \int_{\mathbb R^d} |f_n(y)| |\rho_m(x - y)| dy \overset{\text{H\"older's inequality}}{\le} \|f_n\|_p \|\rho_m\|_{p'}
\end{equation*}
and similarly equicontinuity follows from
\begin{equation*}
\left| \partial_{x_j} h_{n,m}(x) \right| \le \int_{\mathbb R^d} |f_n(y)| |\partial_{x_j} \rho_m(x - y)| dy \overset{\text{H\"older's inequality}}{\le} \|f_n\|_p \|\partial_{x_j} \rho_m\|_{p'}.
\end{equation*}
Furthermore, for $m, n \in \mathbb N$ we have $\operatorname{supp} h_{n,m} \subseteq \overline{S + B_1(0)}$ and $S$ is bounded, which is why our modified Arzel\`a--Ascoli theorem applies to each $\mathcal H_m$ and also to all subsets of $\mathcal H_m$. Now define a function
\begin{equation*}
k: \mathbb N \times \mathbb N \to \mathbb N
\end{equation*}
as thus: $k(n,1)$ is such that $(h_{k(n,1),1})_{n \in \mathbb N}$ is the convergent subsequence of $\mathcal H_1$ as given by our modified Arzel\`a--Ascoli theorem, and if $k(n,m-1)$ is already defined, $k(n,m)$ is such that $(h_{k(n,m),m})_{n \in \mathbb N}$ is the convergent subsequence of $(h_{k(n,m-1)})_{n \in \mathbb N} \subseteq \mathcal H_m$ as given by our modified Arzel\`a--Ascoli theorem. From this, we define a subsequence $(g_n)_{n \in \mathbb N}$ of $(f_n)_{n \in \mathbb N}$ as thus:
\begin{equation*}
g_n := f_{k(n,n)}.
\end{equation*}
We claim that $(g_n)_{n \in \mathbb N}$ is a Cauchy sequence in $L^p(S)$. Indeed, let $\epsilon > 0$. If we set
\begin{equation*}
C_1 := \left( \int_{B_1(0)} |\rho_1(y)|^{p'} dy \right)^{p/p'} ~~~ \text{and} ~~~ C_2 := \int_{B_1(0)} 1 dy
\end{equation*}
where $p'$ is the H\"older conjugate of $p$, we get for all $g \in \mathcal G$
\begin{align*}
\int_S |g * \rho_n(x) - g(x)|^p dx & = \int_S \left| \int_{B_{1/n}(0)} \rho_n(y) g(x - y) dy - \int_{B_{1/n}(0)} \rho_n(y) g(x) dy \right|^p dx \\
& \le \int_S \left( \int_{B_1(0)} |\rho_1(y)| |g(x - y/n) - g(x)| dy \right)^p dx \\
& \overset{\text{H\"older's inequality}}{\le} C_1 \int_S \int_{B_1(0)} |g(x - y/n) - g(x)|^p dy dx \\
& \overset{\text{Fubini's theorem}}{=} C_1 \int_{B_1(0)} \int_S |g(x - y/n) - g(x)|^p dx dy \\
& \le C_1 C_2 \sup_{\|y\| < 1} \int_S |g(x - y/n) - g(x)|^p dx
\end{align*}
and hence, taking the supremum over $g \in \mathcal G$ of this, first on the right hand side and then on the left hand side, we get
\begin{equation*}
\sup_{g \in \mathcal G} \int_S |g * \rho_n(x) - g(x)|^p dx \le C_1 C_2  \sup_{g \in \mathcal G}\sup_{\|y\| < 1} \int_S |g(x - y/n) - g(x)|^p dx.
\end{equation*}
Therefore, by our assumption on countable subsets of $\mathcal F$, we may choose $J \in \mathbb N$ sufficiently large so that for all $g \in \mathcal G$
\begin{equation*}
\|g * \rho_J - g\|_p \le \epsilon/3.
\end{equation*}
Furthermore, by construction, the sequence $(h_{k(J,n),J})_{n \in \mathbb N}$ is a Cauchy sequence in $\mathcal C(\overline{S + B_1(0)})$ and hence also in $L^p(S)$, since $S$ is bounded. Hence we may pick $M \in \mathbb N$ such that for all $n, m \ge M$ we have
\begin{equation*}
\|h_{k(J,n),J} - h_{k(J,m),J}\|_p < \epsilon/3.
\end{equation*}
Then set $N := \max \{J, M\}$ to obtain for $m,n \ge N$ that
\begin{equation*}
\|g_n - g_m\|_p \le \|g_n - h_{k(n,n),J}\|_p + \|h_{k(n,n),J} - h_{k(m,m),J}\|_p + \|h_{k(m,m),J} - g_m\|_p < \epsilon/3 + \epsilon/3 + \epsilon/3 = \epsilon.
\end{equation*}
\end{proof}

\subsection{The relationship to the axiom of choice}
\label{subsec:arzela-ascoli_and_frechet-kolmogorov_relation_to_ac}

We first note that the cardinality of separable metric spaces is always less than or equal to the cardinality of the real numbers $\mathbb R$ (every separable metric space is homeomorphic to a subspace of the Hilbert cube, see for instance Bourbaki \cite[Proposition~12, p.~156]{bourbakiGT2}). Then we note that for a compact $K \subset \mathbb R^d$, the space $\mathcal C(K)$ of continuous, real valued functions with domain of definition $K$ is separable (this may be proven without the axiom of choice by approximating any function in such a space by a suitable multi-dimensional Bernstein polynomial (see for instance \cite[Section~5.2, p.~119ff.]{reimer}) and then in turn approximating the Bernstein polynomial by a rational polynomial, where the set of rational polynomials is countable). Further we note that for any measurable set $S \subseteq \mathbb R^d$, the space $L^p(S)$ may be regarded as a subspace of $L^p(\mathbb R^d)$ by identifying equivalence classes from $L^p(S)$ with equivalence classes of $L^p(\mathbb R^d)$ which are almost everywhere zero on $\mathbb R^d \setminus S$ through the obvious bijective function, which preserves metric. Now the space $L^p(\mathbb R^d)$ is separable (see for instance Brezis \cite[Theorem~4.13, p.~98~f.]{brezis}). Thus, we conclude that the cardinality of both $\mathcal C(K)$ and $L^p(S)$ is less than or equal to the cardinality of the real numbers $\mathbb R$. Furthermore, the functions $\mathbb R \ni x \mapsto x \mathbf 1_K \in \mathcal C(K)$ ($K \neq \emptyset$) and $\mathbb R \ni x \mapsto \left[x \mathbf 1_S\right] \in L^p(S)$ ($S \subseteq \mathbb R^d$ with nonzero measure, square brackets indicating equivalence class formation) are injections, which is why the cardinalities of $\mathcal C(K)$ ($K \neq \emptyset$) and $L^p(S)$ ($S \subseteq \mathbb R^d$ with nonzero measure) are equal to the cardinality of $\mathbb R$ by the Schr\"oder--Bernstein theorem (for the statement and a choiceless proof of the Schr\"oder--Bernstein theorem see for instance Halmos \cite[Chapter~22, p.~88~f.]{halmos}). With this in mind, we now explicate the relationship between the axiom of choice and the Arzel\`a--Ascoli and Fr\'echet--Kolmogorov theorems.

Our investigation uses the same method deployed in Rhineghost \cite{rhineghost} and is thus based on the following result proved by Herrlich and Strecker \cite[Main~theorem, p.~553]{herrlichstrecker}:

\begin{theorem}\label{countable_real_unbounded}
The axiom of countable choice for subsets of real numbers is equivalent to the statement \enquote*{Every unbounded subset of $\mathbb R$ contains a countable, unbounded subset.}
\end{theorem}

Furthermore, we use the following result (see for instance \cite[Remark on p.~290]{husek}):

\begin{theorem}\label{equicontinuous_extension}
Given a bounded, equicontinuous set of functions $\mathcal F \subseteq \mathcal C(K)$ where $K \subset \mathbb R^d$ is compact, one may extend each function in $\mathcal F$ so that a bounded, equicontinuous set of functions defined on $\overline{B_R(0)}$ arises, where $R>0$ is such that $K \subseteq \overline{B_R(0)}$.
\end{theorem}

\begin{theorem}\label{ascoli_equiv_choice_reals}
Assume the validity of ZF. Then the Arzel\`a--Ascoli theorem (theorem \ref{arzela_ascoli}) is true if and only if the axiom of countable choice for subsets of the real numbers $\mathbb R$ is true.
\end{theorem}

\begin{proof}
We first prove the Arzel\`a--Ascoli theorem from the axiom of countable choice. Indeed, $(2) \Rightarrow (1)$ in theorem \ref{arzela_ascoli} is covered by theorem \ref{modified_ascoli}. For $(1) \Rightarrow (2)$, we argue by contradiction: When $\mathcal F$ is not bounded, we use the axiom of countable choice for subsets of the real numbers to obtain an unbounded countable set $\mathcal G \subseteq \mathcal F$ ; when $\mathcal F$ is not equicontinuous, there exists an $\epsilon > 0$ and an $x \in K$ such that
\begin{equation*}
S_n := \left\{ f \in \mathcal F \middle| \exists y \in B_{1/n}(x): |f(x) - f(y)| \ge \epsilon \right\}
\end{equation*}
is nonempty for every $n \in \mathbb N$, and we apply the axiom of countable choice for subsets of the reals to get a countable, non-equicontinuous $\mathcal G \subset \mathcal F$. In both cases, first we extend everything to a suitable $\overline{B_R(0)}$ using theorem \ref{equicontinuous_extension} to ensure separability of $K$, and then apply theorem \ref{modified_ascoli} to obtain a contradiction.

Then we deduce the axiom of countable choice from the Arzel\`a--Ascoli theorem \ref{arzela_ascoli}. Indeed, let $S \subseteq \mathbb R$ be an unbounded subset of $\mathbb R$, and set $K := \overline{B_R(0)}$, where $R > 0$. Set
\begin{equation*}
\mathcal F := \left\{ x \mapsto s \middle| s \in S \right\} \subset \mathcal C(K),
\end{equation*}
the constant functions for the elements of $S$. By the Arzel\`a--Ascoli theorem, $\mathcal F$ contains a sequence which does not have a convergent subsequence, and hence by our modified Arzel\`a--Ascoli theorem \ref{modified_ascoli}, there exists a countable $\mathcal G \subseteq \mathcal F$ that is either not bounded or not equicontinuous. But $\mathcal F$ is equicontinuous, hence $\mathcal G$ is unbounded. The theorem follows from theorem \ref{countable_real_unbounded}.
\end{proof}

\begin{theorem}
Assume the validity of ZF. Then the Fr\'echet--Kolmogorov theorem is true if and only if the axiom of countable choice for subsets of the real numbers $\mathbb R$ is true.
\end{theorem}

\begin{proof}
We first prove the Fr\'echet--Kolmogorov theorem from the axiom of countable choice for subsets of the real numbers. Indeed, $(2) \Rightarrow (1)$ is covered by theorem \ref{modified_frechet}, and for $(1) \Rightarrow (2)$ we suppose for a contradiction that either $\mathcal F$ is unbounded or that $\mathcal F$ does not satisfy
\begin{equation}\label{limit_equiintegral}
\lim_{h \to 0} \sup_{f \in \mathcal F} \int_S |f(x +h) - f(x)|^p dx = 0; \tag{a}
\end{equation}
if $\mathcal F$ is unbounded, we may pick an unbounded countable subset $\mathcal G \subseteq \mathcal F$ (apply countable choice for subsets of the reals to $\{f \in \mathcal F| \|f\|_p > n\}$), and if $\mathcal F$ does not satisfy equation \eqref{limit_equiintegral}, then there exists $\epsilon > 0$ such that for all $n$, the set
\begin{equation*}
S_n := \left\{ f \in \mathcal F \middle| \exists h \in B_{1/n}(0): \int_S |f(x +h) - f(x)|^p dx > \epsilon \right\}
\end{equation*}
is nonempty, and by choosing from these sets, we get a countable $\mathcal G \subseteq \mathcal F$ that does not satisfy \eqref{limit_equiintegral} with $f \in \mathcal G$, and in both cases, theorem \ref{modified_frechet} gives a contradiction.

To prove that the axiom of countable choice for subsets of real numbers follows from the Fr\'echet--Kolmogorov theorem, we proceed exactly as in the proof of theorem \ref{ascoli_equiv_choice_reals}.
\end{proof}

\section{The uniform boundedness principle}
\label{sec:ubp}

\subsection{Introduction}
\label{subsec:ubp_introduction}

The uniform boundedness principle may be stated as thus (see for instance Brezis \cite[Theorem~2.2, p.~32]{brezis}):
\begin{theorem}[Uniform boundedness principle]
Let $(X, \| \cdot \|_X)$ be a Banach space, $(Y, \| \cdot \|_Y)$ a normed space. Let $(T_\alpha)_{\alpha \in A}$ be a family of linear and continuous functions from $X$ to $Y$. If
\begin{equation*}
\forall x \in X: \sup_{\alpha \in A} \|T_\alpha(x)\|_Y < \infty,
\end{equation*}
then
\begin{equation*}
\sup_{\alpha \in A} \|T_\alpha\|_{op} < \infty.
\end{equation*}
\end{theorem}
Roughly speaking, this theorem could be described as asserting: \enquote*{If a family of linear and continuous functions is pointwise bounded, it is also bounded with regard to the operator norm.}

In what follows, we will prove that the axiom of countable choice implies the uniform boundedness principle and that several axioms follow from the uniform boundedness principle. In fact, this will also clarify the choice-axiomatic strength of several other theorems. This is due to the following:

\begin{theorem}\label{ubp_equivalents}
The following are equivalent in ZF:
\begin{enumerate}
\item Every Banach space is barrelled
\item Every lower semi-continuous seminorm on a Banach space is continuous
\item The uniform boundedness principle holds
\end{enumerate}
\end{theorem}

\begin{proof}
\begin{description}
\item [$(1) \Leftrightarrow (2)$] Schechter \cite[27.32 and 27.33, p.~737]{schechter}
\item [$(2) \Rightarrow (3)$] Bourbaki \cite[Theorem~1, p.~III.25]{bourbakiTVS} and Bourbaki \cite[Theorem~4, p.~362]{bourbakiGT1}
\item [$(3) \Rightarrow (1)$] Schechter \cite[27.35, p.~738f.]{schechter}
\end{description}
\end{proof}

\begin{theorem}\label{closed_graph_equivalents}
The following are equivalent in ZF:
\begin{enumerate}
\item The closed graph theorem (i. e., if $X, Y$ are Banach spaces and $T: X \to Y$ is a linear function such that \begin{equation*}
\operatorname{graph} T := \{(x, T(x))|x \in X\} \subset X \times Y
\end{equation*} is closed, then $T$ is continuous) holds.
\item A sequential version of the closed graph theorem (i. e., if $X, Y$ are Banach spaces and $T: X \to Y$ is a linear function such that \begin{equation*}
\operatorname{graph} T := \{(x, T(x))|x \in X\} \subset X \times Y
\end{equation*} is sequentially closed\footnote{Note that the statement \enquote*{In every metric space, a set is closed if and only it is sequentially closed} is equivalent to the axiom of countable choice (see \cite[Theorem~2.1, p.~146]{gutierres}); however, closed sets are always sequentially closed, and for graphs the equivalence holds in ZF, as the theorem shows. In particular, for linear functions between Banach spaces, continuity and sequential continuity are equivalent in ZF.}, then $T$ is continuous) holds.
\item The open mapping theorem (i.e. whenever $X$ and $Y$ are Banach spaces and $T: X \to Y$ is a linear, continuous and surjective function, then $T$ is open) holds.
\end{enumerate}
\end{theorem}

\begin{proof}
\begin{description}
\item [$(3) \Rightarrow (2)$] Brezis \cite[Theorem~2.9, p.~37 and corollaries~7~and~8, p.~35]{brezis}
\item [$(1) \Rightarrow (3)$] Robertson and Robertson \cite[Theorem~3, p.~12]{robertson}\footnote{Note that when $\tau: X \times Y \to X \times Y, \tau(x,y) := (x, y - Tx)$, then $\operatorname{graph} T = \tau^{-1}(X \times \{0\})$.}
\end{description}
\end{proof}

Further, from Schechter \cite[27.34, p.~737f.]{schechter} it follows that any of the statements listed in \ref{closed_graph_equivalents} implies the statements in \ref{ubp_equivalents}.

Our knowledge apart from this article of the relation between the UBP and the axiom of choice stems from an article by Norbert Brunner \cite{brunner}. Indeed, he proved that if every Banach space is barrelled, then the axiom of countable multiple choice (see for instance Herrlich \cite[Definition~2.10, p.~14]{herrlich} or indeed axiom \ref{cmc}) holds (see Brunner \cite[Lemma~4, p.~124f.]{brunner}). Furthermore, he proved that given ZF and the axiom of countable finite choice (i.e. from a sequence of finite sets one can extract a sequence of members), the axiom of countable multiple choice suffices to prove that every Banach space is barrelled (see Brunner \cite[proof of Lemma~5, p.~125f.]{brunner}).

Combining this with theorem \ref{ubp_equivalents} and the fact that trivially, the axiom of countable choice implies both the axiom of countable multiple choice and the axiom of countable finite choice, it follows that
\begin{enumerate}
\item the axiom of countable choice implies the uniform boundedness principle, and
\item the uniform boundedness principle implies the axiom of countable multiple choice.
\end{enumerate}

In what follows, we will give new, direct proofs of these two facts, and further deduce two additional choice-like axioms from the uniform boundedness principle.

\subsection{Countable choice implies UBP}
\label{subsec:ubp_implied_by_cc}

In slightly modifying an argument given by Alan D. Sokal \cite{sokal}, we are able to prove the uniform boundedness principle using nothing more than the Zermelo--Fraenkel system and the axiom of countable choice (as noted above, this has been done before).

\begin{theorem}\label{cc_implies_ubp}
Assume the Zermelo--Fraenkel system and the axiom of countable choice. Then the uniform boundedness principle holds.
\end{theorem}

\begin{proof}
Assume for a contradiction that $(T_\alpha)_{\alpha \in A}$ is an unbounded family of linear, continuous functions from a Banach space $(X, \| \cdot \|_X)$ to a Banach space $(Y, \| \cdot \|_Y)$. Then all the sets
\begin{equation*}
A_n := \left\{ T_\alpha \middle| \alpha \in A, \|T_\alpha\|_{op} > 4^n \right\}
\end{equation*}
are nonempty. Hence, by the axiom of countable choice, we may pick a sequence $(T_n)_{n \in \mathbb N}$ such that for each $n \in \mathbb N$ we have $T_n \in A_n$. By definition of the operator norm, all the sets
\begin{equation*}
B_n := \left\{ x \in X \middle| \|x\|_X \le 1, \|T_n(x)\|_Y > \frac{2}{3} \|T_n\|_{op} \right\}
\end{equation*}
are nonempty. A second application of the axiom of countable choice hence permits us to choose a sequence $(x_n)_{n \in \mathbb N}$ such that for all $n \in \mathbb N$, we have $x_n \in B_n$.

Now define a function by
\begin{equation*}
f: \mathbb N \times X \to X,~~ f(n, x) := \begin{cases}
x + 3^{-(n+1)} x_{n+1} & \left\|T_{n+1}\left(x + 3^{-(n+1)} x_{n+1}\right)\right\|_Y > 3^{-(n+1)} \frac{2}{3} \|T_{n+1}\|_{op} \\
x - 3^{-(n+1)} x_{n+1} & \text{else}.
\end{cases}
\end{equation*}

We claim that for each $n \in \mathbb N$ there exists exactly one $n$-tuple $(z_{n,1}, \ldots, z_{n,n})$ such that
\begin{enumerate}
\item $z_{n,1} = x_1$
\item $z_{n,k+1} = f(k,z_{n,k})$ for $k \in \{1, \ldots, n-1\}$.
\end{enumerate}
Existence is proved by induction on $n$; for uniqueness, assume otherwise and use that $f$ is a function and hence can have only one value. For a given $n$, we have by induction on $k$ that $z_{n,k} = z_{n-1,k}$ for $k \in \{1, \ldots, n-1\}$. Define a sequence $(y_n)_{n \in \mathbb N}$ by $y_n := z_{n,n}$. We get for $n \in \mathbb N$
\begin{equation*}
y_{n+1} = z_{n+1,n+1} = f(n, z_{n+1,n}) = f(n, z_{n,n}) = f(n, y_n).
\end{equation*}
The sequence $(y_n)_{n \in \mathbb N}$ has two properties:

\begin{enumerate}
\item When $k, n \in \mathbb N$, we have
\begin{equation*}
\|y_n - y_{n+k}\| \le \sum_{j=0}^\infty \|y_{n+j+1} - y_{n+j}\| \le \sum_{j=0}^\infty 3^{-(n+j+1)} = 3^{-(n+1)} \frac{3}{2} = 3^{-n}\frac{1}{2}.
\end{equation*}
\item For $x \in X$ (since the maximum is larger than the average and due to the triangle inequality)
\begin{align*}
\max\left\{ \left\| T_n(x + 3^{-n} x_n) \right\|, \left\| T_n(x - 3^{-n}  x_n) \right\| \right\} & \ge \frac{1}{2} \left( \left\| T_n(x + 3^{-n} x_n) \right\| + \left\| T_n(x - 3^{-n}  x_n) \right\| \right) \\
& \ge 3^{-n} \|T_n(x_n)\| \ge 3^{-n} \frac{2}{3} \|T_n\|_{op}
\end{align*}
and hence
\begin{equation*}
\|T_n(y_n)\| \ge 3^{-n} \frac{2}{3} \|T_n\|_{op}.
\end{equation*}
\end{enumerate}

From the first property, $(y_n)_{n \in \mathbb N}$ is a Cauchy sequence, hence convergent to some $y \in X$. Then for $k, n \in \mathbb N$
\begin{equation*}
\|y_n - y\| \le \|y_n - y_k\| + \|y_k - y\|,
\end{equation*}
and letting $k \to \infty$ proves, together with the first property, that $\|y_n - y\| \le 3^{-n}/2$.

Combining this with the second property, we get
\begin{equation*}
\|T_n(y)\| \ge \|T_n(y_n)\| - \|T_n(y_n - y)\| \ge 3^{-n} \frac{2}{3} \|T_n\|_{op} - 3^{-n}\frac{1}{2} \|T_n\|_{op} = \frac{1}{6} 3^{-n} \|T_n\|_{op} > \frac{1}{6} (4/3)^n \to \infty;
\end{equation*}
that is, $(T_\alpha)_{\alpha \in A}$ is unbounded in $y$.
\end{proof}

Note that instead of using the axiom of countable choice in selecting the functions $(T_n)_{n \in \mathbb N}$, we could have instead defined the $B_n$ as
\begin{equation*}
B_n := \left\{ x \in X \middle| \|x\|_X \le 1, \exists \alpha \in A: \|T_\alpha(x)\|_Y > \frac{2}{3} \|T_\alpha\|_{op} > \frac{2}{3} 4^n \right\}
\end{equation*}
and only used the axiom of countable choice on the $B_n$. Then we would have replaced all the inequalities by inequalities for which there exists a suitable $\alpha$ such that they hold.

\subsection{Axioms implied by the UBP}
\label{subsec:ubp_implied_axioms}

\subsubsection{General proof strategy}

Clearly, as became apparent in the proof of the last subsection, the uniform boundedness principle may be reformulated as thus:

\begin{theorem}[Uniform boundedness principle, reformulated]
Let $(T_\alpha)_{\alpha \in A}$ be an unbounded family of linear functions from a Banach space $X$ to a normed space $Y$. Then there exists $x \in X$ such that $\{ T_\alpha(x) | \alpha \in A \} \subseteq Y$ is unbounded.
\end{theorem}

Hence, we see that given an unbounded family of linear, continuous functions, the uniform boundedness principle translates into an \emph{existence} claim, namely the existence of a point where the respective family of linear, continuous functions is unbounded. If we therefore are able to associate to such a point subsets of a given family of sets, we can use an existence result as such to obtain variants of the axiom of choice. In fact, using the right Banach spaces, this will be possible. The spaces considered in this paper are created using the following construction:

Assume we are given a countable family of Banach spaces $X_1, X_2, \ldots, X_n, \ldots$. Then we may construct from them new Banach spaces in which $X_1, X_2, \ldots, X_n, \ldots$ are isometrically contained in a canonical fashion. Namely, if $1 \le p \le \infty$, set
\begin{equation*}
\sideset{}{^p}\bigoplus_{n \in \mathbb N} X_n := \left\{ (x_n)_{n \in \mathbb N} \in \prod_{n \in \mathbb N} X_n \middle| \left\| (x_n)_{n \in \mathbb N} \right\|_p < \infty \right\}
\end{equation*}
where the norm $\| \cdot \|_p$ is given by
\begin{equation*}
\left\| (x_n)_{n \in \mathbb N} \right\|_p := \begin{cases}
\left( \sum\limits_{n=1}^\infty \|x_n\|_{X_n}^p \right)^{1/p} & p < \infty \\
\sup\limits_{n \in \mathbb N} \|x_n\|_{X_n} & p = \infty
\end{cases},
\end{equation*}
where for each $n$ $\| \cdot \|_{X_n}$ is the norm of $X_n$. All these spaces are Banach spaces (see for instance Helemskii \cite[p. 127]{helemskii}). The space
\begin{equation*}
\sideset{}{^p}\bigoplus_{n \in \mathbb N} X_n
\end{equation*}
is called the \emph{$\ell^p$ sum} of $X_1, X_2, \ldots, X_n, \ldots$.

In what follows, we will associate to a sequence of sets $(S_n)_{n \in \mathbb N}$ spaces $(X_n)_{n \in \mathbb N}$. The choice axioms that we want to deduce from the uniform boundedness principle will assert that given a sequence $(S_n)_{n \in \mathbb N}$, for infinitely many $n$ we can choose sets $M_n \subseteq S_n$ whose cardinality obeys a certain restriction (for instance is finite or obeys some bound). In order to execute the deduction, we pick the spaces $(X_n)_{n \in \mathbb N}$ such that, for a given $n$, all elements of $X_n$ having a certain property (for instance, nonzero elements) yield a set $M_n$ obeying the desired condition; this is the case when the elements of $X_n$ with the certain property are all sufficiently \enquote*{asymmetric} in their structure. Then we use the uniform boundedness principle, applied to a suitable unbounded family of linear, continuous functions (which will be defined on a suitable $\ell^p$ sum of the $(X_n)_{n \in \mathbb N}$), to get elements of $X_n$ for infinitely many $n$ which have the certain property that will yield a suitable $M_n$.

In our exposition, we will dedicate a subsubsection to each axiom that will be deduced from the UBP. In each subsubsection, we will start with a lemma explaining why the elements with the \emph{certain} property are sufficiently asymmetric such that subsets of the desired cardinality may be selected.

\subsubsection{Choosing finite subsets}
\label{subsubsec:ubp_implies_cmc}

We will now present another proof for the fact that the uniform boundedness principle implies the axiom of countable multiple choice, which is defined as thus (see for instance \cite[Definition~2.10, p.~14]{herrlich}):

\begin{axiom}[Countable multiple choice]\label{cmc}
Let $(S_n)_{n \in \mathbb N}$ be a sequence of sets. Then there exists a sequence of nonempty sets $(M_n)_{n \in \mathbb N}$ such that for all $n$ we have $M_n \subseteq S_n$ and $M_n$ is finite.
\end{axiom}

We will use the fact that this axiom is equivalent to the following seemingly modifieder axiom:

\begin{axiom}[Partial countable multiple choice]
Let $(S_n)_{n \in \mathbb N}$ be a sequence of sets. Then there exists an infinite $I \subseteq \mathbb N$ and a family of nonempty sets $(M_n)_{n \in I}$ such that for all $n \in I$ we have $M_n \subseteq S_n$ and $M_n$ is finite.
\end{axiom}

The proof of equivalence of the two axioms is essentially the same as the proof of equivalence of the axiom of partial countable choice to the axiom of countable choice (see for instance \cite[Theorem~2.12~3., p.~15]{herrlich}) and follows easily from Keremedis \cite[Lemma~1.2, p.~570]{keremedis}.

Let $S$ be an arbitrary set. $\mathcal P(S)$, the power set of $S$, is a $\sigma$-algebra, and the \emph{counting measure} on $S$ is defined as
\begin{equation*}
\mu(E) = \begin{cases}
\# E & \# E < \infty \\
\infty & \text{otherwise}
\end{cases},
\end{equation*}
where $E \subseteq S$ is arbitrary (see for instance Tao \cite[Example 1.4.26, p. 90f.]{taoMeasure}). $L^1(S, \mu)$ is a Banach space (see for instance Brezis \cite[Theorem 4.8, p. 93]{brezis}).

We note the following lemma regarding $L^1(S, \mu)$:

\begin{lemma}\label{nonzero_to_finite}
Let $(S_\alpha)_{\alpha \in A}$ be a family of sets, and let $\mu_\alpha$ be the counting measure on $S_\alpha$. There exists a function
\begin{equation*}
\Phi_{(S_\alpha)_{\alpha \in A}}: \bigcup_{\alpha \in A} \left\{ f \in L^1(S_\alpha, \mu_\alpha) \middle| \int_{S_\alpha} f(\sigma) d\mu_\alpha(\sigma) \neq 0 \right\} \to \bigcup_{\alpha \in A} \left\{ T \subseteq S_\alpha \middle| 0 < \#T < \infty \right\}
\end{equation*}
such that whenever $f \in L^1(S_\alpha, \mu_\alpha)$, then $\Phi_{(S_\alpha)_{\alpha \in A}}(f) \subseteq S_\alpha$.
\end{lemma}

\begin{proof}
Let $\alpha \in A$ and $f \in L^1(S_\alpha, \mu_\alpha)$ such that $\int_{S_\alpha} f(\sigma) d\mu_\alpha(\sigma) \neq 0$ be given. We partition the split real number line $\mathbb R \setminus \{0\}$ into countably many subsets as thus:
\begin{equation*}
\mathbb R \setminus \{0\} = \bigcup_{n \in \mathbb Z} \left( \left[ -2^{n+1} , -2^n \right) \cup \left( 2^n , 2^{n+1} \right] \right).
\end{equation*}
The following three observations are immediate:
\begin{enumerate}
\item For all $n \in \mathbb Z$, at most finitely many elements of $f(S_\alpha)$ are in $\left[ -2^{n+1} , -2^n \right) \cup \left( 2^n , 2^{n+1} \right]$.
\item $\left[ -2^{n+1} , -2^n \right) \cup \left( 2^n , 2^{n+1} \right]$ contains a point of $f(S_\alpha)$ only for finitely many \emph{positive} $n$.
\item $\mathbb R \setminus \{0\}$, and hence at least one of the sets $\left[ -2^{n+1} , -2^n \right) \cup \left( 2^n , 2^{n+1} \right]$, must contain a point of $f(S_\alpha)$.
\end{enumerate}
Hence, we choose $n \in \mathbb Z$ maximal such that there are some points of $f(S_\alpha)$ in $\left[ -2^{n+1} , -2^n \right) \cup \left( 2^n , 2^{n+1} \right]$; these are finitely many, and we define $\Phi_{(S_\alpha)_{\alpha \in A}}(f)$ to be the set of these points.
\end{proof}

\begin{theorem}\label{ubp_implies_pcmc}
Assume that the uniform boundedness principle holds. Then, given a sequence $(S_n)_{n \in \mathbb N}$ of sets, there exists a sequence of sets $(M_n)_{n \in \mathbb N}$ such that for all $n \in \mathbb N$ $M_n \subseteq S_n$ and $|M_n| < \infty$ and for infinitely many $n$, $M_n \neq \emptyset$.
\end{theorem}

\begin{proof}
For each $n \in \mathbb N$, set $X_n := L^1(S_n,\mu_n)$, where $\mu_n$ is the counting measure on $S_n$. Then set
\begin{equation*}
X := \sideset{}{^p}\bigoplus_{n \in \mathbb N} X_n,
\end{equation*}
where the choice of $1 \le p \le \infty$ does not matter; we may for instance take $p=1$. Set $Y := \mathbb R$, where the norm is given by the absolute value of the reals. For each $n$ define a linear function $T_n: X \to Y$ by
\begin{equation*}
T_n \left( (x_n)_{n \in \mathbb N} \right) := 4^n \int_{S_n} x_n(\sigma) d \mu_n(\sigma).
\end{equation*}
If, for a fixed $n$, a $\sigma \in S_n$ is picked, and the element $x = (x_k)_{k \in \mathbb N} \in X$ defined by the sequence
\begin{equation*}
x_k = \begin{cases}
0 & n \neq k \\
\mathbf 1_{\{\sigma\}} & n = k
\end{cases}
\end{equation*}
is considered (where $\mathbf 1_A$ is defined to be the indicator function on the set $A$), it becomes evident that $\|T_n\|_{op} \ge 4^n$; that is, the family of linear functions $(T_n)_{n \in \mathbb N}$ is uniformly unbounded. Hence, by the uniform boundedness principle, it is also pointwise unbounded. In particular, we find $x = (x_n)_{n \in \mathbb N} \in X$ such that
\begin{equation*}
\int_{S_n} x_n(\sigma) d \mu_n(\sigma) \neq 0
\end{equation*}
for infinitely many $n$, say for all $n \in I$, where $I \subseteq \mathbb N$ is infinite. For all such $n$, we then define $M_n := \Phi_{(S_n)_{n \in I}}(x_n)$, where $\Phi_{(S_n)_{n \in I}}$ is as in lemma \ref{nonzero_to_finite}.
\end{proof}

Hence, the uniform boundedness principle implies the axiom of partial countable multiple choice, and thus the axiom of countable multiple choice.

\subsubsection{Choosing subsets of asymptotically bounded cardinality}
\label{subsubsec:ubp_implies_some_asymptotic_choice}

Now we will deduce from the uniform boundedness principle the following choice axiom:

\begin{axiom}[Axiom of partial countable asymptotic choice]\label{asymptotic_axiom}
Let $(S_n)_{n \in \mathbb N}$ be a sequence of finite sets, and let $(\lambda_n)_{n \in \mathbb N}$ be a sequence of positive real numbers which is unbounded. Then there exists an infinite set $I \subseteq \mathbb N$, a family of sets $(M_n)_{n \in I}$ and a constant $C > 0$ such that for all $n \in I$ we have $M_n \subseteq S_n$, $M_n \neq \emptyset$ and $|M_n| \le C \lambda_n$.
\end{axiom}

The deduction of this axiom is based on the following lemma:

\begin{lemma}\label{inequality_cardinality_bound}
Let $S$ be a finite set, $\mu$ the counting measure on $S$, and $f: S \to \mathbb R$ a real-valued function on $S$. If for a $C > 0$ we have
\begin{equation*}
\int_S |f(\sigma)| d\mu(\sigma) \le C \sup_{\sigma \in S} |f(\sigma)|,
\end{equation*}
then $\# \left\{ \sigma \in S \middle| |f(\sigma)| = \sup_{\sigma \in S} |f(\sigma)| \right\} \le C$.
\end{lemma}

\begin{proof}
Assume otherwise. Then
\begin{equation*}
\int_S |f(\sigma)| d\mu(\sigma) \ge \#\left\{ \sigma \in S \middle| |f(\sigma)| = \sup_{\sigma \in S} |f(\sigma)| \right\} \cdot \sup_{\sigma \in S} |f(\sigma)| > C \sup_{\sigma \in S} |f(\sigma)|.
\end{equation*}
\end{proof}

\begin{theorem}
The uniform boundedness principle implies axiom \ref{asymptotic_axiom}.
\end{theorem}

\begin{proof}
Let $(S_n)_{n \in \mathbb N}$ be a sequence of finite, nonempty sets, and let $(\lambda_n)_{n \in \mathbb N}$ be an unbounded sequence of positive, real numbers. Let $\mu_n$ be the counting measure on $S_n$. Set $X_n := L^1(S_n, \mu_n)$ and $Y_n := L^\infty(S_n, \mu_n)$ and then
\begin{equation*}
X := \sideset{}{^\infty}\bigoplus_{n \in \mathbb N} X_n ~~~\text{and}~~~ Y := \sideset{}{^p}\bigoplus_{n \in \mathbb N} Y_n;
\end{equation*}
the choice of $1 \le p \le \infty$ doesn't matter; e.g. $p=1$. Define a sequence of linear functions by
\begin{equation*}
T_n: X \to Y, T_n \left( (x_k)_{k \in \mathbb N} \right) := \lambda_n \left( \delta_{n k} x_k \right)_{k \in \mathbb N},
\end{equation*}
where
\begin{equation*}
\delta_{n k} = \begin{cases}
1 & n = k \\
0 & \text{else}
\end{cases}
\end{equation*}
is the Kronecker delta. If, for a fixed $n$, a $\sigma \in S_n$ is picked, and the element $x = (x_k)_{k \in \mathbb N} \in X$ defined by the sequence
\begin{equation*}
x_k = \begin{cases}
0 & n \neq k \\
\mathbf 1_{\{\sigma\}} & n = k
\end{cases}
\end{equation*}
is considered, it becomes evident that $\|T_n\|_{op} \ge \lambda_n$. It follows that the family of linear functions $(T_n)_{n \in \mathbb N}$ is uniformly unbounded, and hence, by the uniform boundedness principle, pointwise unbounded. Hence pick $(x_k)_{k \in \mathbb N} \in X$ such that the set $\left\{ \left\|T_n\left( (x_k)_{k \in \mathbb N} \right) \right\|_Y \middle| n \in \mathbb N\right\}$ is an unbounded subset of the real numbers. In particular, there exists an infinite set $I \subseteq \mathbb N$ such that for all $n \in I$ we have
\begin{equation*}
1 \le \left\|T_n\left( (x_k)_{k \in \mathbb N} \right) \right\|_Y = \lambda_n \left\| \left( \delta_{n k} x_k \right)_{k \in \mathbb N} \right\|_Y = \lambda_n \max_{\sigma \in S_n} |x_n(\sigma)|
\end{equation*}
But since $(x_k)_{k \in \mathbb N} \in X$, there exists a $C > 0$ such that
\begin{equation*}
\forall n \in \mathbb N: C \ge \|x_n\|_{X_n} = \int_{S_n} |x_n(\sigma)| d\mu_n(\sigma).
\end{equation*}
Hence, for all $n \in I$ we have
\begin{equation*}
\int_{S_n} |x_n(\sigma)| d\mu_n(\sigma) \le C \le C \lambda_n \max_{\sigma \in S_n} |x_n(\sigma)|.
\end{equation*}
By lemma \ref{inequality_cardinality_bound}, by defining $M_n := \left\{ \sigma \in S_n \middle| |x_n(\sigma)| = \sup_{\sigma \in S_n} |x_n(\sigma)| \right\}$ for all $n \in I$ we get a family of sets as required by axiom \ref{asymptotic_axiom}.
\end{proof}

\subsubsection{Choosing singletons}
\label{subsubsec:ubp_implies_bounded_choice}

\begin{axiom}[Axiom of countable choice for sets of bounded, finite cardinality]\label{axiom_bounded_choice}
The axiom of countable choice for sets of bounded, finite cardinality, or for short $\mathbf{CC}(\{1, \ldots, n\})$, shall mean the following: If $(S_k)_{k \in \mathbb N}$ is a sequence of sets such that $\forall k \in \mathbb N: \# S_k \le n$, then
\begin{equation*}
\prod_{n \in \mathbb N} S_k \neq 0.
\end{equation*}
\end{axiom}

\begin{definition}\label{simplex_spaces_def}
Let $S$ be a set of finite cardinality, say $\# S = n \in \mathbb N$. Then define
\begin{align*}
U_S &:= \left\{ (x_\sigma)_{\sigma \in S} \middle| \forall \sigma \in S: x_\sigma \in \mathbb R \right\} \\
V_S &:= \{t(1)_{\sigma \in S} | t \in \mathbb R\} \\
W_S &:= U_S / V_S.
\end{align*}
\end{definition}

Trivially, the following holds:

\begin{lemma}\label{subset_from_nonzero}
Let $S$ be a set of finite cardinality and $w = (x_\sigma)_{\sigma \in S} + V_S \in W_S$ such that $w \neq 0$. Then the set
\begin{equation*}
\left\{ \sigma \in S \middle| x_\sigma \neq \min_{\sigma \in S} x_\sigma \right\}
\end{equation*}
is independent of the representative of $w$ and a nonempty, proper subset of $S$.
\end{lemma}

\begin{theorem}
Assume that the uniform boundedness principle is true. Then axiom \ref{axiom_bounded_choice} holds.
\end{theorem}

\begin{proof}
Let a family of sets $(S_k)_{k \in \mathbb N}$ be given such that for all $k$ we have $\# S_k \le n$. For each $k$, we form the spaces $U_{S_k}$ and $W_{S_k}$ as given in definition \ref{simplex_spaces_def}. These are finite-dimensional real vector spaces and hence may be normed to obtain Banach spaces. We set $X_k$ to be the Banach space that results from norming $U_{S_k}$, and $Y_k$ the Banach space that results from norming $W_{S_k}$. Then we define spaces
\begin{equation*}
X := \sideset{}{^p}\bigoplus_{n \in \mathbb N} X_n, ~~~ Y := \sideset{}{^q}\bigoplus_{n \in \mathbb N} Y_n
\end{equation*}
(where the choice of $1 \le p, q \le \infty$ doesn't matter; we may as well take $p = q = 1$) and linear functions
\begin{equation*}
T_k: X \to Y, T_k\left((x_m)_{m \in \mathbb N}\right) := 4^k \left( \delta_{k,m} \pi_{W_{S_m}}(x_m) \right)_{m \in \mathbb N},
\end{equation*}
where
\begin{equation*}
\pi_{W_{S_m}}: U_{S_m} \to W_{S_m}
\end{equation*}
are the canonical projections. Once more, considering the element $x = (x_m)_{m \in \mathbb N} \in X$ defined by the sequence
\begin{equation*}
x_m = \begin{cases}
0 & k \neq m \\
\mathbf 1_{\{\sigma\}} & k = m
\end{cases}
\end{equation*}
yields that the family of linear functions $(T_k)_{k \in \mathbb N}$ is unbounded, an we take the uniform boundedness principle to get a point $(x_m)_{m \in \mathbb N}$ such that
\begin{equation*}
\pi_{W_{S_m}}(x_m) \neq 0 \text{ infinitely often.}
\end{equation*}
By lemma \ref{subset_from_nonzero} this gives an infinite $I \subseteq \mathbb N$ and a family $(M_n)_{n \in I}$ of nonempty sets such that for all $n$ we have $M_n \subsetneq S_n$. Repeating this process $n$ times yields the theorem.
\end{proof}

\subsection{The UBP for real-valued functions}
\label{subsec:weak_ubp_equiv_cmc}

There is a weak form of the uniform boundedness principle that is equivalent (in the Zermelo--Fraenkel system) to the axiom of countable multiple choice. Which is:

\begin{theorem}[Weak form of the uniform boundedness principle]\label{weak_ubp}
Let $X$ be a Banach space, and let $(T_n)_{n \in \mathbb N}$ be a sequence of linear and continuous functions from $X$ to $\mathbb R$ such that for all $x \in X$, the set
\begin{equation*}
\left\{T_k(x) \middle| k \in \mathbb N \right\}
\end{equation*}
is a bounded subset of $\mathbb R$. Then
\begin{equation*}
\sup_{k \in \mathbb N} \|T_k\|_{op} < \infty.
\end{equation*}
\end{theorem}

\begin{theorem}
Assuming the Zermelo--Fraenkel system, but not any choice, theorem \ref{weak_ubp} is equivalent to the axiom of countable multiple choice.
\end{theorem}

\begin{proof}
For necessity, we note that the proof of theorem \ref{ubp_implies_pcmc} goes through without modification. For sufficiency, note that in our situation, we may modify the proof of theorem \ref{cc_implies_ubp} as follows so that only countable multiple choice is needed:
\begin{enumerate}
\item The first use of the axiom of countable choice is avoided by taking the first element of the sequence $(T_k)_{k \in \mathbb N}$ which obeys the desired bound.
\item The axiom of countable multiple choice is applied to get a finite subset $S_n \subset B_n$ ($B_n$ as in the proof of theorem \ref{cc_implies_ubp}) for each $n$, and if we set
\begin{equation*}
x_n := \frac{1}{|S_n|} \sum_{x \in S_n} x,
\end{equation*}
then $\|x_n\| \le 1$ and $T_n(x_n) > \frac{2}{3} \|T_n\|_{op}$.
\end{enumerate}
\end{proof}

\subsection{Outlook}
\label{subsec:ubp_outlook}

We have derived several choice-like axioms from the uniform boundedness principle. I believe that in conjunction, they are not strong enough to prove the axiom of countable choice, in particular in view of the fact that even if for each $n \in \mathbb N$ and each family $(S_\alpha)_{\alpha \in A}$ of sets of cardinality $n$
\begin{equation*}
\prod_{\alpha \in A} S_\alpha \neq 0,
\end{equation*}
the axiom of finite choice $\mathbf{AC}(\operatorname{fin})$ (see for instance Herrlich \cite[Definition~2.6, p.~14]{jech}) does not follow (see for instance Jech \cite[Theorem~7.11, p.~107]{jech}). I have been unable to deduce the axiom of countable choice from the uniform boundedness principle, but there are several lines of attack that seem promising:
\begin{enumerate}
\item The use of a modified form of ultraproducts.
\item The use of different characterisations of the uniform boundedness principle (for instance, as mentioned it is equivalent in ZF to the statement that all Banach spaces are barrelled).
\item Using the fact that $Y$ merely needs to be a normed space.
\end{enumerate}
Now the gap between what we proved and the axiom of countable choice is less huge than one would perhaps suspect; many choice axioms that we deduced from the UBP are of the \enquote*{partial type}, but since partial countable choice is equivalent to countable choice, it would suffice to prove that from a sequence of sets $(S_n)_{n \in \mathbb N}$ that has \emph{arbitrary asymptotic behaviour}, one can pick a sequence of subsets $(M_n)_{n \in \mathbb N}$ such that for infinitely many $n$ we have $\# M_n < C$ for a $C > 0$ and $M_n \neq \emptyset$.

\section{Acknowledgements}

I would like to thank two users of a mathematical Q\&A website who helped me by explaining to me some facts about cardinalities. Furthermore, I would like to thank one of the referees for providing me with the wonderful and simple algebraic description of the spaces of definition \ref{simplex_spaces_def}, which I had before described inefficiently (and not wholly correctly) by analytical means, and both editors for all the valuable suggestions they made. I would also like to thank the editor for several explanations, and the Journal of Logic and Analysis for providing a venue without charge for papers of this kind.

\bibliographystyle{plain}
\bibliography{bib-paper}{}

\begin{thebibliography}{10}

\bibitem{abbott}
Stephen Abbott.
\newblock {\em Understanding Analysis}, volume 207 of {\em Undergraduate Texts
  in Mathematics}.
\newblock Springer New York, 2 edition.

\bibitem{bentleyHerrlich}
H.~L. Bentley and H.~Herrlich.
\newblock Countable choice and pseudometric spaces.
\newblock {\em Topology and its Applications}, 85(1-3):153--164, 5 1985.

\bibitem{bourbakiGT2}
Nicolas Bourbaki.
\newblock {\em General Topology 2}, volume~3 of {\em Elements of Mathematics}.
\newblock Addison-Wesley, 1966.

\bibitem{bourbakiGT1}
Nicolas Bourbaki.
\newblock {\em General Topology 1}, volume~3 of {\em Elements of Mathematics}.
\newblock Springer-Verlag Berlin Heidelberg, 3 edition, 1995.

\bibitem{bourbakiTVS}
Nicolas Bourbaki.
\newblock {\em Topological Vector Spaces}, volume~5 of {\em Elements of
  Mathematics}.
\newblock Springer Berlin Heidelberg, 2003.

\bibitem{brezis}
Haim Brezis.
\newblock {\em Functional Analysis, Sobolev Spaces and Partial Differential
  Equations}.
\newblock Universitext. Springer New York, 1 edition, 2011.

\bibitem{brunner}
Norbert Brunner.
\newblock Garnir's dream spaces with hamel bases.
\newblock {\em Archiv f\"ur mathematische Logik und Grundlagenforschung},
  26(1):123--126, 12 1987.

\bibitem{gutierres}
Gon{\c c}alo Gutierres.
\newblock On countable choice and sequential spaces.
\newblock {\em Mathematical Logic Quarterly}, 54(2):145--152, 4 2008.

\bibitem{halmos}
Paul~R. Halmos.
\newblock {\em Naive Set Theory}.
\newblock Undergraduate Texts in Mathematics. Springer New York, 1974.

\bibitem{helemskii}
A.~Ya. Helemskii.
\newblock {\em Lectures and Exercises on Functional Analysis}, volume 233 of
  {\em Translations of Mathematical Monographs}.
\newblock American Mathematical Society, 2006.

\bibitem{herrlichAscoli}
Horst Herrlich.
\newblock
  \href{http://ftp.math.uni-rostock.de/pub/romako/heft51/herrlich.ps.Z}{The
  Ascoli Theorem is equivalent to the Boolean Prime Ideal Theorem}.
\newblock {\em Rostocker Mathematisches Kolloquium}, 51:137--140, 1997.

\bibitem{herrlich}
Horst Herrlich.
\newblock {\em Axiom of Choice}.
\newblock Lecture Notes in Mathematics. Springer Berlin Heidelberg, 2006.

\bibitem{herrlichstrecker}
Horst Herrlich and George~E. Strecker.
\newblock
  \href{http://dml.cz/bitstream/handle/10338.dmlcz/118952/CommentatMathUnivCarolRetro_38-1997-3_11.pdf}{When
  is $\mathbb{N}$ Lindel\"of?}
\newblock {\em Commentationes Mathematicae Universitatis Carolinae},
  38(3):553--556, 1997.

\bibitem{husek}
Miroslav Hu{\v s}ek.
\newblock
  \href{http://www.eweb.unex.es/eweb/extracta/Vol-25-3/25J3Husek.pdf}{Extension
  of Mappings and Pseudometrics}.
\newblock {\em Extracta Mathematicae}, 25(3):277--308, 12 2010.

\bibitem{jech}
Thomas~J. Jech.
\newblock {\em The Axiom of Choice}, volume~75 of {\em Studies in Logic and the
  Foundations of Mathematics}, chapter~7, pages 97--118.
\newblock North Holland, 1973.

\bibitem{keremedis}
Kyriakos Keremedis.
\newblock Disasters in topology without the axiom of choice.
\newblock {\em Arch. Math. Logic}, 40(8):569--580, 11 2001.

\bibitem{reimer}
Manfred Reimer.
\newblock {\em Multivariate Polynomial Approximation}, volume 144 of {\em
  International Series of Numerical Mathematics}.
\newblock Springer Basel, 2003.

\bibitem{rhineghost}
Y.~T. Rhineghost.
\newblock The naturals are lindel{\"o}f iff ascoli holds.
\newblock In J{\"u}rgen Koslowski and Austin Melton, editors, {\em Categorical
  Perspectives}, Trends in Mathematics, chapter~11, pages 191--194.
  Birkh\"auser Boston, 2001.

\bibitem{robertson}
Alexander~P. Robertson and Wendy Robertson.
\newblock On the closed graph theorem.
\newblock {\em Glasgow Mathematical Journal}, 3(1):9--12, 12 1956.

\bibitem{rudin}
Walter Rudin.
\newblock {\em Principles of Mathematical Analysis}.
\newblock International Series in Pure \& Applied Mathematics. McGraw-Hill, 3
  edition, 1976.

\bibitem{schechter}
Eric Schechter.
\newblock {\em Handbook of Analysis and Its Foundations}, chapter~27, pages
  721--751.
\newblock Elsevier, 1997.

\bibitem{sokal}
Alan~D. Sokal.
\newblock A really simple elementary proof of the uniform boundedness theorem.
\newblock {\em The American Mathematical Monthly}, 118(5):450--452, 5 2011.

\bibitem{taoReal}
Terence Tao.
\newblock {\em
  \href{https://terrytao.files.wordpress.com/2012/12/gsm-117-tao3-epsilon1.pdf}{An
  Epsilon of Room, I: Real Analysis}}.
\newblock Number 117 in Graduate Studies in Mathematics. American Mathematical
  Society, 2010.

\bibitem{taoMeasure}
Terence Tao.
\newblock {\em
  \href{https://terrytao.files.wordpress.com/2012/12/gsm-126-tao5-measure-book.pdf}{An
  introduction to measure theory}}.
\newblock Number 126 in Graduate Studies in Mathematics. American Mathematical
  Society, 2011.

\bibitem{yosida}
Kôsaku Yosida.
\newblock {\em Functional Analysis}.
\newblock Classics in Mathematics. Springer Berlin Heidelberg, 6 edition, 1995.

\end{thebibliography}

\end{document}